\documentclass[12pt,a4paper]{article}

\usepackage{indentfirst, amsfonts, amsmath, amsthm, amssymb, amscd}
\usepackage{graphics, latexsym}
\usepackage{epsfig}
\usepackage{latexsym}

\newtheorem{theorem}{\bf Theorem}[section]

\newtheorem{proposition}[theorem]{\bf Proposition}
\newtheorem{remark}[theorem]{\bf Remark}

\font\bbbld=msbm10 scaled\magstep1
\newcommand{\bfR}{\hbox{\bbbld R}}
\newcommand{\bfC}{\hbox{\bbbld C}}

\newcommand{\bfH}{\hbox{\bbbld H}}

\newcommand{\bfP}{\hbox{\bbbld P}}

\newcommand{\bfS}{\hbox{\bbbld S}}

\def\G{{\Gamma}}
\def\la{\langle}
\def\ra{\rangle}

\def\G{{\Gamma}}

\def\Sc{{\cal S}}

\def\Nb{{\bf N}}

\def\Xb{{\bf X}}

\def\Xb{{\bf X}}

\def\G{\Gamma}

\def\o{{\omega}}

\def\ts{\tilde{S}}

\def\sl{\bfP SL(2,\bfC)}

\begin{document}

\title{CMC-1 surfaces in hyperbolic $3$-space using the
Bianchi-Cal\`{o} method}
\author{Levi Lopes de Lima\thanks{Departamento de Matem\'atica, Universidade Federal do Cear\'a, Fortaleza, Brasil.}\,\,\,\, Pedro Roitman\thanks{Departamento de Matem\'atica, Universidade Federal do Cear\'a, Fortaleza, Brasil.}}

\date{\today}

\maketitle

\begin{abstract}

In this  note we present a method for constructing CMC-1 surfaces in hyperbolic $3$-space
$\bfH^3(-1)$ in terms of holomorphic data first introduced in Bianchi's Lezioni di Geometria Differenziale of 1927, therefore predating  by many years the modern approaches due to Bryant, Small and others. Besides its obvious historical interest, this note aims to complement Bianchi's analysis by deriving explicit formulae for CMC-1 surfaces and comparing the various approaches encountered in the literature.

\end{abstract}

\section{Introduction}
\label{intr}

It is generally accepted that the theory of surfaces in hyperbolic 3-space $\bfH^3(-1)$ with
constant mean curvature equal to one (CMC-1 surfaces, for short) started with
a seminal paper by R. Bryant (\cite{b}), where he derives  a representation for such
surfaces in terms of holomorphic data. More precisely, Bryant's recipe works as follows. Start with a holomorphic curve $\Sc\subset\sl$ which is null with respect to the conformal structure inherited from the Killing-Cartan form and recall the natural projection $\pi:\sl\to\bfH^3(-1)$, 
$\pi(\omega)=\omega\overline{\omega}^t$, where $\bfH^3(-1)$ is realized via the hermitian model and $\sl$ is the group of orientation preserving  isometries of $\bfH^3(-1)$. It turns out  that, at least locally, any CMC-1 surface is of the form $\pi(\Sc)$ for some such $\Sc$. 
In this setting, the primary object is a pair of holomorphic functions defining another (holomorphic) curve on the Lie algebra of $\sl$ from which $\Sigma$ is obtained after solving a certain first order differential equation. Thus  an integration procedure appears in the description and  
this partially accounted for the inherent difficulties in constructing examples at the very  early stages of the theory. But notice however the striking analogy with the well-known Weierstrass
representation for minimal surfaces in $\bfR^3$ (\cite{l}), as these are obtained as real slices of null curves in $\bfC^3$.
 
After the appearance of Bryant's investigation, many other researchers  contributed to the subject. For example, starting with the paper \cite{uy}, M. Umehara and K. Yamada refined substantially Bryant's approach and were able to construct a varied class of examples of CMC-1 surfaces, besides developing many interesting global aspects in the theory. On the other hand, 
inspired by Hitchin's generalization (\cite{h}) of   
Penrose's twistor theory, A. J. Small (\cite{s}) gave a clean characterization of null curves in $\sl$. In fact, Small realizes $\sl$ as the complement of the quadric $Q_2=\{ad-bc=0 \}\subset\bfP^3$ and observes that twistors lead naturally to the construction of a Gauss transform $\Gamma_{\Sc}$ of a  null curve $\Sc\subset\sl$ as a curve lying in $Q_2^*$, the quadric dual to $Q_2$. Moreover, he shows that, aside from the technical issue of missing points lying in $Q_2$, one has 
\begin{equation}\label{small}
\Sc=\Gamma_{\Sc}^*,
\end{equation}
where $\Gamma_{\Sc}^*$ is the dual to $\Gamma_{\Sc}$.  This gives a classical   
procedure for recovering $\Sc$ from $\Gamma_{\Sc}$ and in this context one is naturally compelled to regard  the Gauss transform  as the primary object. The final blow is the well-known fact that  $Q_2^*=Q_2=\bfP\times\bfP$, a product of two projective lines,  and we can locally   represent 
\begin{equation}\label{sma}
\Gamma_{\Sc}=(f,g),
\end{equation}
for holomorphic functions $f$ and $g$ on the Riemann surface $M$ underlying $\G_{\Sc}$. Since dualization is carried out by applying algebraic operations on the derivatives up to second order of $(f,g)$  one ends up with an explicit formula for CMC-1 surfaces involving no integration whatsoever!   
We stress that this approach has the obvious drawback of being local in nature but the elimination of quadratures reminds us of an alternative (and much lesser known) way of describing the classical Weiertrass representation, which was rediscovered in modern times by Hitchin (\cite{h2}) and Small (\cite{s2}) again by using twistorial methods.

The above historical account reflects an attempt to interpret the prevalent view in the CMC-1 surfaces community. Consensual as it may be, the main purpose of  this note is to point out that it is not entirely correct. In effect, in Bianchi's Lezioni di Geometria Differenziale (\cite{b}), edited in 1927, we may find a recipe
for constructing CMC-1 surfaces out of holomorphic data.  
However, the motivation for writing this note goes   beyond this historical curiosity,
for it seems that there are at least two other reasons for exhibiting this 
old method to a wider audience. 

The first one is that the method  allows
one to start with an arbitrary holomorphic map 
$
f 
$
defined in a region  $\Omega\subset \bfC$ and, elaborating upon Bianchi's ideas, to end up with explicit formulae for a
CMC-1 surface (see Theorem \ref{cent2}). Moreover, this map $f$ has
an immediate geometric interpretation: it is simply the parametrized
hyperbolic Gauss map, or in other words, the expression for the hyperbolic
Gauss map in terms of a local complex parameter on $\Omega$. 
In addition, the resulting formulae involve algebraic expressions in terms of the derivatives of $f$ up to second order  and no integration at all. As such, they come close to the spirit of Small's  approach as described above. In fact, and this is an important issue here, our formulae coincide with Small's  if the pertinent transformation between models for $\bfH^3(-1)$ is carried out (see Section \ref{bcmeth}).

The second reason is that in his way toward the   construction of CMC-1
surfaces, Bianchi translates to hyperbolic geometry the solution of  a strictly Euclidean-geometric problem  involving the
rolling of a pair of isometric surfaces, thereby establishing a surprising  linking between these two geometries. As a bonus, his approach is entirely elementary in the sense that for its understanding it is only required familiarity with the fundamentals of three definitely classical disciplines: Euclidean geometry, hyperbolic geometry and complex function theory. This should be compared to the heavy algebraic-geometric machinery used by Small.

This article is organized as follows.
In Section \ref{clasec} we review some classical concepts in order to
describe the rolling problem in Euclidean geometry. In particular we derive, following Bianchi, the so called Cal\`{o}'s formulae. In Section \ref{bcmeth} we relate the solution of the Euclidean problem of Section \ref{clasec} to the construction of CMC-1 surfaces in $\bfH^3(-1)$. We obtain an explicit parametrization for a CMC-1 surface in terms of a conformal map. Finally in Section \ref{examp} we exhibit some simple examples constructed via the
Bianchi-Cal\`{o} method.

\section{Classical concepts and the Bianchi-Cal\`{o} method}
\label{clasec}

\subsection{Congruence of spheres and envelopes.}
\label{conse}

Here is the first classical concept we shall meet. A {\it congruence of spheres} is a smooth two-parameter family of spheres
in $\Bbb{R}^3,$ that we will suppose parametrized by coordinates $(u,v)$. To each such congruence we may associate a function $R=R(u,v)$, the {\it radius function},  describing the radii of the spheres  in the congruence. We also assume that the vector function $\Xb=\Xb(u,v)$ describing the centers of the spheres defines a regular surface which we call the
{\it surface of centers}. Such a congruence of spheres will be denoted by $[\Xb,R]$.

Generically there are two surfaces, the so-called {\it envelopes}, associated to a given congruence $[\Xb,R]$. In effect,  a point  $p\in\bfR^3$ belongs to an envelope $\xi $ if $p\in S$ for some sphere $S$ in $[\Xb,R]$ and moreover $T_p\Sigma =T_pS$. The next proposition gives the
expression for the
envelopes 
in terms of  the unit normal vector $\Nb=\Nb(u,v)$ and the metric of 
the surface of centers.

\begin{proposition}
In coordinates $(u,v)$, 
\begin{equation}
\xi =\Xb-R\left(\Delta(\Xb,R)\pm \sqrt{1-\Delta
_1R}\,\Nb\right),  \label{form}
\end{equation}
where 
\begin{equation}
\Delta (\Xb,R)=(R_uA_{11}+R_vA_{12})\Xb_u+(R_uA_{21}+R_vA_{22})
\Xb_v,  \label{firm}
\end{equation}
and 
\begin{equation}
\Delta _1R=R_u^2A_{11}+2R_uR_vA_{12}+R_v^2A_{22}.  \label{farm}
\end{equation}
Here, the matrix $A=[A_{ij}]$ is the inverse of the matrix defined by the
metric in the given coordinates.
\end{proposition}

\begin{proof} The conditions defining the envelopes are
expressed as 
\begin{equation}
{\xi }={\Xb}-R{\nu },  \label{env}
\end{equation}
\begin{equation}
\la\nu ,{\xi }_u\ra=\la{\nu },{\xi }_v\ra=0,
\label{derienv}
\end{equation}
where ${\nu }$ is the unit vector in the direction of  
$
{\Xb}-{\xi }$, so that taking derivatives of (\ref{env}) and using the fact that $\left| {\nu }\right|
=1$ we obtain 
\begin{equation}
\la{\nu },{\Xb}_u\ra=R_u, \,\,\, \la{\nu },{\Xb}_v\ra=R_v \label{mix} 
\end{equation}
Now, write 
\begin{equation}\label{expnu}
{\nu }=a{\Xb}_u+b{\Xb}_v+c{\Nb},
\end{equation}
take inner products of this with ${\Xb}_u$ and ${\Xb}_v$ and solve the
linear system 
\begin{eqnarray*}
g_{11}a+g_{12}b=R_u, \\
g_{21}a+g_{22}b=R_v,
\end{eqnarray*}
so as to obtain
\begin{eqnarray*}
a=A_{11}R_u+A_{12}R_v, \\
b=A_{21}R_u+A_{22}R_v.
\end{eqnarray*}
Finally, $c$ is further determined  using that $\left| \nu \right| =1$:
\begin{equation}
c=\pm \sqrt{1-\Delta _1R}.  \label{cc}
\end{equation}
Substitution yields the desired formula for the envelopes.
\end{proof}

If a congruence of spheres has two distinct envelopes we then have a
natural correspondence between their points, namely,  points on
distinct envelopes correspond if they are the contact points of the envelopes with a
given sphere of the congruence.
We will allow the degenerate case where one
of the envelopes reduces to a point and notice  that even in this situation  the unit vector ${\nu }$
appearing in the above proposition is still well defined, a fact we shall use in the following 
result due originally to Beltrami.

\begin{proposition}\label{bel}
Let $[\Xb,R]$ be a congruence of
spheres with  ${\xi }$ being 
one of its envelopes (possibly degenerated to a point), obtained say by choosing the positive sign in (\ref{form}), and consider the angles 
$\omega _1,\omega _2$ and $ \sigma $ formed between the unit vector ${\nu }$ in the direction of 
${\Xb}-{\xi }$ and ${\Xb}_u$, ${\Xb}_v$ and ${\Nb}$, respectively.
Suppose further that another surface 
$\tilde{{\Xb}}$ isometric to ${\Xb}$ is given and  consider the  congruence
of spheres $[\tilde{{\Xb}},R]$. Finally, let $\tilde{\xi }$ be the envelope for $[\tilde{\Xb},R]$ obtained by choosing the same sign as we did in order to get ${\xi }$ starting from $[\Xb,R]$.  Then the  angles between the unit vector $\tilde{\nu }$ in the
direction of $\tilde{\Xb}-\tilde{\xi }$ and $\tilde{\Xb}_u,
\tilde{\Xb} _{v}$ and $\tilde{\Nb}$, respectively, where $\tilde{\Nb}$ is the unit normal to $\tilde{\Xb}$, coincide with the corresponding angles for $[\Xb,R]$.
\end{proposition}

\begin{proof} 
From (\ref{mix}), (\ref{expnu}) and (\ref{cc}) one has
\begin{eqnarray*}
\cos \omega _1 & = & \frac{\la{\nu },{X}_u\ra}{\sqrt{g_{11}}} = 
\frac{R_u}{\sqrt{g_{11}}}, \\
\cos \omega _2 & = & \frac{\la{\nu },{X}_v\ra}{\sqrt{g_{22}}}  = 
\frac{R_v}{\sqrt{g_{22}}}, \\
\cos \sigma & = & \la{\nu },{N}\ra  =  \sqrt{1-\Delta _1R},
\end{eqnarray*}
and since the surfaces of centers of the congruences, which have the same radius function, are isometric to each other, the right hand sides above are the same when we consider the angles of both congruences.
\end{proof}

\medskip
\centerline{\epsfxsize=15cm \epsfysize=8.4cm \epsfbox{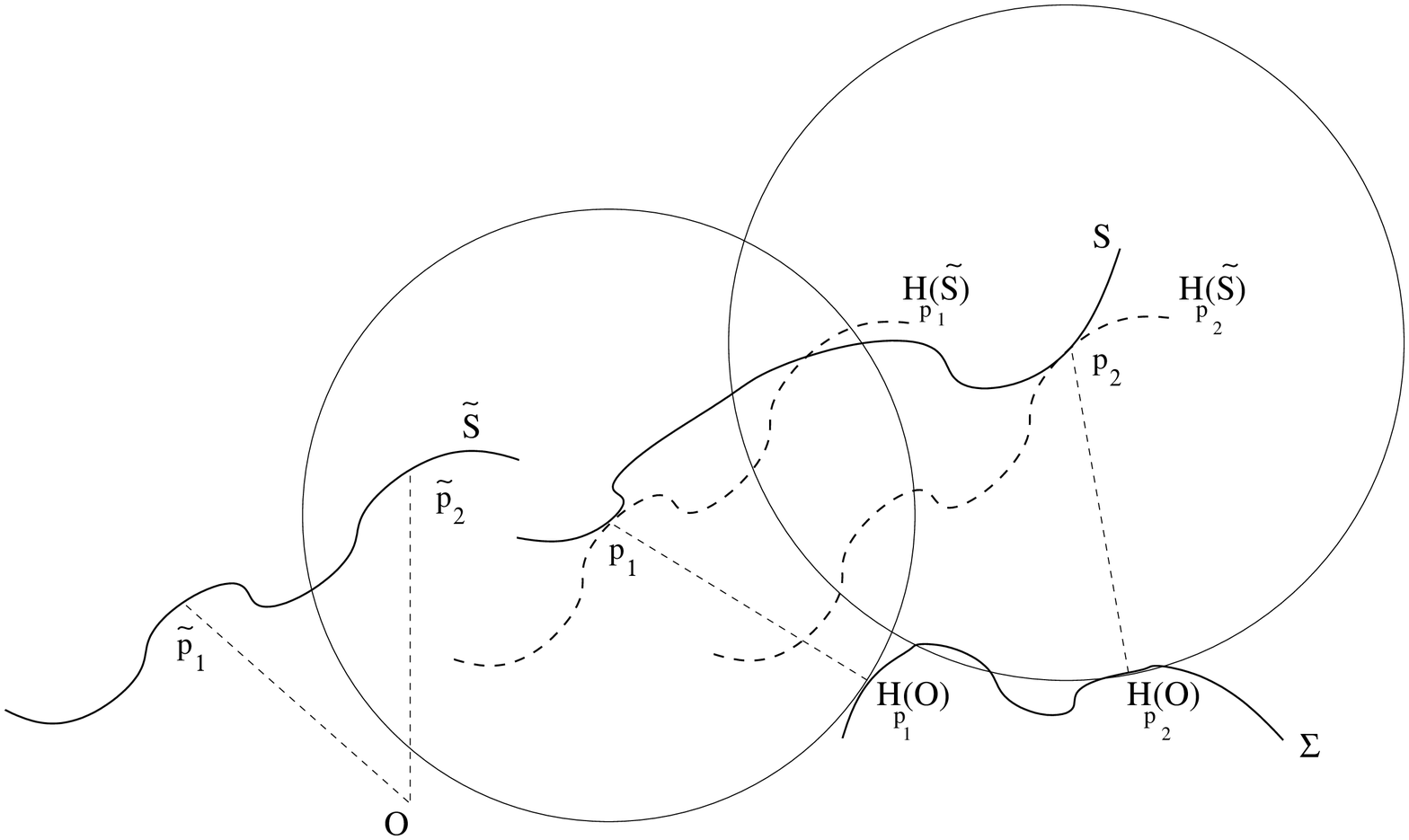}}   
\smallskip
\centerline{Fig. 2.1}
\medskip
\subsection{Rolling of isometric surfaces}
\label{roll}

We now describe the rolling of isometric surfaces and show how a congruence
of spheres is associated to such a rolling. Consider a pair $(S,\ts)$ of
isometric surfaces in $\Bbb{R}^3,$ and let $p$ $\in S$ and $\tilde{p}\in \ts$ be points
corresponding under the isometry. Suppose $S$ is fixed in space and
consider the  two-parameter family of positions of congruent copies of $\ts$
such that to each $p\in S$ we consider a rigid motion of $\Bbb{R}^3$ (call it
$H_p$) sending 
$\tilde{p}$ to $p$, $T_{\tilde{p}}\ts$ to $T_pS$, and further adjusted so that the differential of the
isometry composed with $H_p$ is the identity map. This two-parameter
family of positions for copies of $\ts$ is called the {\it rolling }of 
$\ts$ over $S$. The surface $\ts$ is called the {\it rolled surface} and 
$S$ the {\it support surface}.

Now fix a point $O$ $\in \Bbb{R}^3$ and consider its image under the  two-parameter
family of rigid motions associated to the rolling of $\ts$ over $S$. In
the generic case, the motion of $O$ defines a surface $\Sigma$ called the  {\it  rolling surface} with respect to the {\it satellite point} $O$.

The crucial point now is that
the two concepts introduced so far, namely congruence of spheres and rolling of surfaces, share a close relationship. 
More precisely, we have

\begin{proposition}
Given a   rolling of $\ts$ over $S$ as above,  the rolling 
surface $\Sigma $ can be viewed as an envelope of a congruence of spheres
having $S$ as its surface of centers and the sizes of
the corresponding line segments joining points of $\ts$ to $O$ as radii.
\end{proposition}

\begin{proof}
Look at $O$ as a degenerate envelope
corresponding to the congruence of spheres having $\ts$ as surface of centers
and passing through $O$, see 
figure 2.1. By Proposition \ref{bel}, an envelope of the congruence with $S$ as the surface of centers, and same radius
function as the congruence just considered, has the property that  the unit
vector joining a point of it to the corresponding point of $S$ makes the same
angles (with respect to the obvious  fixed basis)
than the corresponding unit vector joining a point $p\in \ts$ to $O$ makes with
the corresponding basis. But this
shows that the point of the envelope coincides with the 
point of the rolling surface.
\end{proof}

\subsection{The Cal\`o's formulae}
\label{calooo}

We are now in a position to formulate the problem in Euclidean geometry that
will lead us, according to Bianchi, to a method for constructing CMC-1 surfaces in $\bfH^3(-1)$:  

\medskip

{\it Find pairs} $(S,\ts)$\textit{
\ of isometric surfaces such that, for a convenient satellite point }$O$, \textit{\ the
rolling surface }$\Sigma $\textit{\ is contained in a plane.}

\medskip

To solve this problem we  work with Cartesian coordinates $(x,y,z)$ in $\bfR^3$, suppose that the plane in question is $\{z=0\}$ and moreover that  the satellite point $O$
is the origin of our coordinate system. Let $\tilde{p}=(\tilde{x},\tilde{y},\tilde{z})$
and $p=(x,y,z)$ denote corresponding points for $\ts$ and $S$, respectively. The necessary
and sufficient conditions for a pair $(S,\ts)$ to be a solution to our
problem are 
\begin{eqnarray}
d\tilde{x}^2+d\tilde{y}^2+d\tilde{z}^2& = & dx^2+dy^2+dz^2,  \label{dxxx}\\
  \tilde{x}^2+\tilde{y}^2+\tilde{z}^2 & = & z^2, \label{x2}  
\end{eqnarray}
the first condition expressing that $S$ and $\ts$ are isometric, and the
second one coming  from the assumption that the radius of the sphere of the congruence
equals the distance from $\tilde{p}$ to $O$.

We now write $\ts$ in polar coordinates $(R,\theta ,\phi )$ centered at $O=(0,0,0)$
so that 
\begin{equation}\label{saf}
\tilde{x} =R\sin \theta \cos \phi ,\,\,\,\, \tilde{y} =R\sin \theta \sin \phi,\,\,\,\, 
\tilde{z} =R\cos \theta. 
\end{equation}
Locally, $\ts$ can be parametrized by  an open subset of the unit sphere $\bfS^2$ centered at $O$ via the inverse of central projection so that if
we write the metric of $\ts$ in these coordinates we get
\[
d\tilde{x}^2+d\tilde{y}^2+d\tilde{z}^2=dR^2+R^2(d\theta ^2+\sin ^2\theta d\phi ^2). 
\]
On the other hand, (\ref{x2}) and (\ref{saf}) imply $R^2=z^2$ and choosing
$R=z$, so that  $S$ now lies in the {\it upper} half-space, we obtain, after using (\ref
{dxxx}), 
\begin{equation}
R^2(d\theta ^2+\sin ^2\theta d\phi ^2)=dx^2+dy^2,  \label{metro}
\end{equation}
which we may interpret as follows: {\it the central projection of} $\ts$ {\it over }
$\bfS^2$ {\it and the orthogonal projection of} $S$ {\it onto the
plane} $\{z=0\}$ {\it define a conformal map from} $\bfS^2$ {\it to the plane} $\{z=0\}$. Thus we make our first contact with complex function theory. 

Via stereographic projection  we introduce on $\bfS^2$ the complex variable 
\begin{equation}\label{param}
\tau =\cot \frac{\theta}{ 2}\, e^{i\phi }, 
\end{equation}
so that
in terms of this parameter one gets  the coordinates of $\tilde{S}$:
\begin{equation*}
\tilde{x} =R\frac{\tau +\bar{\tau}}{\left| \tau \right| ^2+1}, \,\,\,\,\,\,
\tilde{y} =\frac Ri\frac{\tau -\bar{\tau}}{\left| \tau \right| ^2+1}, \,\,\,\,\,\,
\tilde{z}=R\frac{\left| \tau \right| ^2-1}{\left| \tau \right| ^2+1}.
\end{equation*}

Much in the same vein, we can parametrize the corresponding piece in $S$ via a local inverse of the orthonormal projection by means of the complex parameter
$
\zeta =x+iy. 
$
The discussion above allows us to consider $\zeta $ as
a holomorphic function of $\tau $: 
$
\zeta =f(\tau ), 
$
 see figure 2.2. In other words, $f$ describes the isometry between the surfaces in terms of the above chosen complex parameters.

\medskip
\centerline{\epsfxsize=15cm \epsfysize=4.6cm \epsfbox{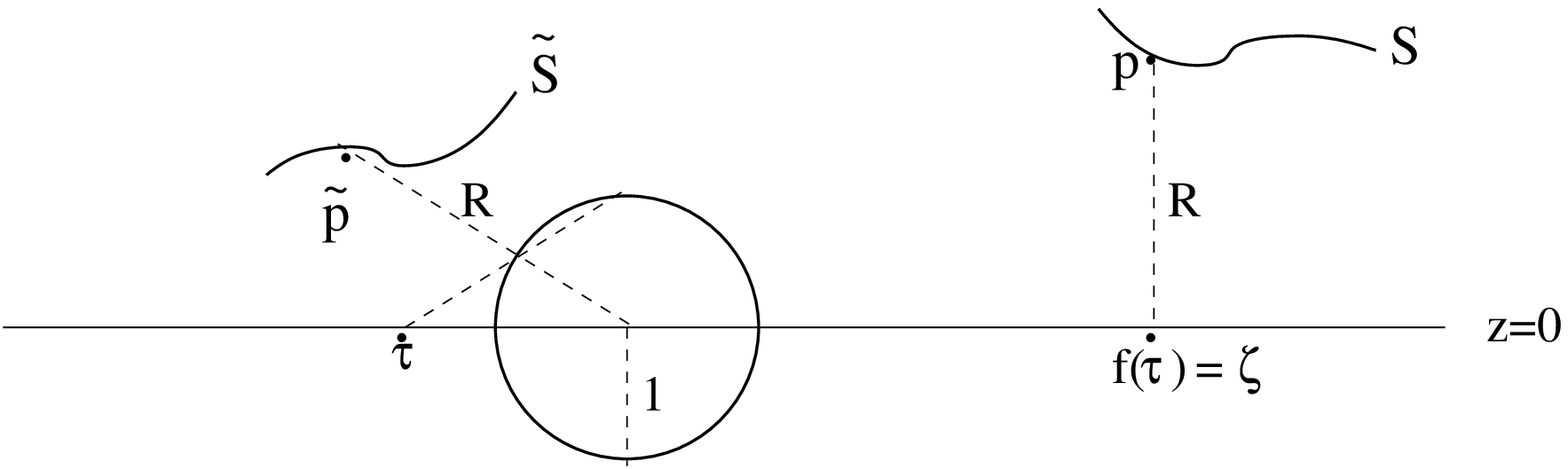}}   
\smallskip
\centerline{Fig. 2.2}
\medskip

The final step is to determine $R$ in terms of $f$. Using (\ref{metro}) and (\ref{param}) we
get 
\[
\frac{4R^2\left| d\tau ^2\right| }{(\left| \tau \right| ^2+1)^2}=\left|
d\zeta ^2\right| =\left| f^{\prime }(\tau )\right| ^2\left| d\tau ^2\right|
, 
\]
from which we conclude that 
\begin{equation}\label{rr}
R=\frac{1+\left| \tau \right| ^2}2\left| f^{\prime }(\tau )\right| , 
\end{equation}
and we have finally met the solution to our problem, namely, the coordinates of $\ts$
and $S$ are respectively given in terms of the holomorphic data as

\begin{equation}\label{calo1}
\tilde{x} = \left| f^{\prime }(\tau )\right| \frac{\tau +\bar{\tau}}2, \,\,\,\,\,\tilde{y} =\left| f^{\prime }(\tau )\right| \frac{\tau -\bar{\tau}}{2i}, \,\,\,\,\,   
\tilde{z} =\left| f^{\prime }(\tau )\right| \frac{\left| \tau \right| ^2-1}2, 
\end{equation}

and

\begin{equation}\label{calo2}
x =\hbox{Re}f(\tau ),\,\,\,\,\, y =\hbox{Im}f(\tau ), \,\,\,\,\,
z =\left| f^{\prime }(\tau )\right| \frac{\left| \tau \right| ^2+1}2.
\end{equation}

The above expressions are called {\it Cal\`{o}'s formulae} since they have been originally 
published by B. Cal\`{o} in 1899 (\cite{c}) in another context involving isometric
surfaces.

\section{The Bianchi-Cal\`{o} method}
\label{bcmeth}

In last section, starting with a holomorphic map $f$, we have determined a
pair of isometric surfaces such that one of the envelopes of the associated congruence of spheres was a plane. Notice that we
also arranged things so that the surface of centers was placed in the upper half-space. Now, in
principle we could also determine the second envelope of the congruence
associated to the rolling, which is then also contained in the upper half space. It can be shown that the correspondence
between the envelopes of the congruence associated to the  Cal\`o's pair $S,\tilde{S}$
considered in the last section is a conformal map. This is proved in
Bianchi's Lezioni when he considers Darboux congruencies and is one of the ingredients in the proof of the following central result, also due to Bianchi.

\begin{theorem}\label{central}
To each pair $(\tilde{S},S)$ of isometric surfaces such that the rolling surface $\Sigma$ of the rolling of $\tilde{S}$ over $S$ is a plane there corresponds a CMC-1 given by the second envelope of the associated  
congruence of spheres considered as a surface in the standard upper half space model of $\bfH^3(-1)$.

\end{theorem}

\begin{proof}
Look at the spheres of
our congruence as horospheres by using  the upper half-space containing the 
surface of centers as a model for $\bfH^3(-1)$. In this way, the correpondence between the envelopes
becomes the hyperbolic Gauss map for the second envelope, which is known to
be conformal, exception made for totally umbilical surfaces, exactly when
the surface is a CMC-1 surface (see \cite{br} or \cite{b}).
\end{proof}

Although Bianchi indicates how one can find CMC-1 surfaces starting with an
arbitrary horlomorphic map $f$ via Theorem \ref{central}, he does not complete his analysis by deriving explicit formulae. In the following, we carry out the calculations and exhibit a
CMC-1 surface in $\bfH^3(-1)$ in terms of the map $f$. The method is indeed very simple: we use Cal\`{o}'s formulae (\ref{rr}) and (\ref{calo2}) to compute the surface  of centers
$S$ and the radius function $R$, then we calculate the envelopes of this congruence of
spheres, one of them being a piece of the plane $\{z=0\}$ and the other one being
our CMC-1 surface. In the end of this section, we check that these formulae coincide with Small' s  and this proves that they define a CMC-1 surface indeed and moreover that any CMC-1 arises locally in this way. We insist however that a proof of these assertions can be given in an elementary way and  entirely avoiding Small's results: one just has to check directly the validity
of the above statements regarding conformality. We stick to the approach presented in the sequel just to stress the equivalence between 
the two methods. 

\begin{theorem}\label{cent2}
In the above situation, the parametrization of a CMC-1 surfaces  in terms of $f$ is given by
\begin{eqnarray}
x &=&{{\rm{Re}}}\,f-\frac{\left| f^{\prime }\right| ^2{{\rm{Re}}}\left(f^{\prime }\tau
\right)+\frac{1+\left| \tau \right| ^2}2{{\rm{Re}}}\left((f^{\prime })^2\bar{f}^{\prime
\prime }\right)}{\left| f^{\prime }\right| ^2+{{\rm{Re}}}\left( f^{\prime }\bar{f}
^{\prime \prime }\bar{\tau}\right) +\frac{\left| f^{\prime \prime }\right|
^2(\left| \tau \right| ^2+1)}4},  \nonumber \\
y &=&{{\rm{Im}}}f-\frac{\left| f^{\prime }\right| ^2{{\rm{Im}}}(f^{\prime }\tau
)+\frac{1+\left| \tau \right| ^2}2{{\rm{Im}}}\left((f^{\prime })^2\bar{f}^{\prime
\prime }\right)}{\left| f^{\prime }\right| ^2+{{\rm{Re}}}\left( f^{\prime }\bar{f}
^{\prime \prime }\bar{\tau}\right) +\frac{\left| f^{\prime \prime }\right|
^2(\left| \tau \right| ^2+1)}4},  \label{bicalo} \\
z &=&\frac{\left| f^{\prime }\right| ^3}{\left| f^{\prime }\right| ^2+{{\rm{
Re}}}\left( f^{\prime }\bar{f}^{\prime \prime }\bar{\tau}\right) +\frac{
\left| f^{\prime \prime }\right| ^2(\left| \tau \right| ^2+1)}4}.  \nonumber
\end{eqnarray}

\end{theorem}

\begin{proof}
In terms of $\tau=u+iv$, $f=f_1+if_2$ and the radius function $R$, the surface $S$ is written as 
\[
\mathbf{X}(\tau )=(f_1,f_2,R), 
\]
so that the coefficients of
the metric become 
$$
g_{11}  = \left| f^{\prime }\right| ^2+R_u^2,\,\,\,
g_{12} = R_uR_v, \,\,\,
g_{22} =\left| f^{\prime }\right| ^2+R_v^2,
$$
where we used the Cauchy-Riemann equations for $f$.

The corresponding determinant is
\[
\chi =g_{11}g_{22}-g_{12}^2=\left| f^{\prime }\right| ^2\left(\left| f^{\prime
}\right| ^2+R_u^2+R_v^2\right)=\left| f^{\prime }\right| ^2\left(\left| f^{\prime
}\right| ^2+\left| \nabla R\right| ^2\right). 
\]
and moreover  
\[
\mathbf{X}_u\wedge \mathbf{X}_v=\left(\alpha _1,\alpha _2,\left| f^{\prime
}\right| ^2\right), 
\]
where 
\begin{equation}\label{alp}
\alpha _1 = R_vf_{2,u}-R_uf_{2,v}, \,\,\,
\alpha _2 = R_uf_{1,v}-R_vf_{1,u},
\end{equation}
so that the unit normal vector is
\[
\mathbf{N}=\frac 1{\sqrt{\chi }}\left(\alpha _1,\alpha _2,\left| f^{\prime
}\right| ^2\right). 
\]

On the other hand, the inverse of the matrix associated to the metric is 
$$
A_{11} = \frac{\left| f^{\prime }\right| ^2+R_v^2}\chi ,\,\,\,
A_{12} = \frac{-R_uR_v}\chi , \,\,\,
A_{22} = \frac{\left| f^{\prime }\right| ^2+R_u^2}\chi .
$$
From (\ref{firm}) and (\ref{farm}) we obtain 
\begin{eqnarray*}
\Delta (X,R) &  = & \frac{R_u\left| f^{\prime }\right| ^2}\chi
\left(f_{1,u},f_{2,u},R_u\right)+\frac{R_v\left| f^{\prime }\right| ^2}\chi
\left(f_{1,v},f_{2,v},R_v\right)  \\
  & = & \frac{\left| f^{\prime }\right| ^2}\chi \left(-\alpha
_1,-\alpha _2,\left| \nabla R\right| ^2\right), 
\end{eqnarray*}
and
\begin{eqnarray*}
\Delta _1R & = & R_u^2\frac{\left| f^{\prime }\right| ^2+R_v^2}\chi +2R_uR_v\frac{
-R_uR_v}\chi +R_v^2\frac{\left| f^{\prime }\right| ^2+R_u^2}\chi \\
& = & \frac{
\left| f^{\prime }\right| ^2\left| \nabla R\right| ^2}\chi , 
\end{eqnarray*}
so that in particular, 
\[
\sqrt{1-\Delta _1R}=\frac{\left| f^{\prime }\right| ^2}{\sqrt{\chi }}. 
\]

In order to calculate the envelopes we use (\ref{form}) in the form
\[
\mathbf{\xi }_{\pm }=\mathbf{X}-R\left(\Delta (\mathbf{X},R)\pm \sqrt{1-\Delta _1R}\,\Nb
\right), 
\]
so that $\mathbf{\xi }_{+}$, the envelope contained in $\{z=0\}$, is given by 
\begin{eqnarray*}
\mathbf{\xi }_{+} & = &\left(f_1,f_2,R\right)-R\left(\frac{\left| f^{\prime }\right| ^2}\chi
\left(-\alpha _1,-\alpha _2,\left| \nabla R\right| ^2\right)+\frac{\left| f^{\prime
}\right| ^2}{\sqrt{\chi }}\frac 1{\sqrt{\chi }}\left(\alpha _1,\alpha _2,\left|
f^{\prime }\right| ^2\right)\right) \\ 
& = & \left(f_1,f_2,R\right)-R\left(\frac{\left| f^{\prime }\right| ^2}\chi
\left(0,0,\left| \nabla \phi \right| ^2+\left| f^{\prime }\right|
^2\right)\right) =  \left(f_1,f_2,0\right), 
\end{eqnarray*}
as expected, and
$\mathbf{\xi }_{-}$, our CMC-1 surface, is 
\begin{eqnarray*}
\mathbf{\xi }_{-} & = & \left(f_1,f_2,R\right)-R\left(\frac{\left| f^{\prime }\right| ^2}\chi
\left(-\alpha _1,-\alpha _2,\left| \nabla R\right| ^2\right)-\frac{\left| f^{\prime
}\right| ^2}{\sqrt{\chi }}\frac 1{\sqrt{\chi }}\left(\alpha _1,\alpha _2,\left|
f^{\prime }\right| ^2\right)\right)\\ 
& = & \left(f_1+\frac{2\alpha _1R}{\left| f^{\prime }\right|
^2+\left| \nabla R\right| ^2},f_2+\frac{2\alpha _2R}{\left| f^{\prime
}\right| ^2+\left| \nabla R\right| ^2},\frac{2R\left| f^{\prime }\right| ^2}{
\left| f^{\prime }\right| ^2+\left| \nabla R\right| ^2}\right). 
\end{eqnarray*}
An explicit formula entirely in terms of $f$ can be found if we use (\ref{rr}) to obtain
\begin{eqnarray*}
R_u &=&u\left| f^{\prime }\right| +\frac{1+\left| \tau\right| ^2}{2\left|
f^{\prime }\right| }{{\rm{Re}}}\left(f^{\prime }\bar{f}^{\prime \prime }\right), \\
R_v &=&v\left| f^{\prime }\right| +\frac{1+\left| \tau\right| ^2}{2\left|
f^{\prime }\right| }{{\rm{Im}}}\left(f^{\prime }\bar{f}^{\prime \prime }\right).
\end{eqnarray*}
It then follows  from (\ref{alp}) that
\begin{eqnarray*}
\alpha _1 &=&\left(v\left| f^{\prime }\right| +\frac{1+\left| \tau\right| ^2}{
2\left| f^{\prime }\right| }{{\rm{Im}}}\left(f^{\prime }\bar{f}^{\prime \prime
}\right)\right)f_{2,u}-\left(u\left| f^{\prime }\right| +\frac{1+\left| \tau\right| ^2}{2\left|
f^{\prime }\right| }{{\rm{Re}}}\left(f^{\prime }\bar{f}^{\prime \prime }\right)\right)f_{1,u} \\
&=&-\left| f^{\prime }\right|{{ \rm{Re}}}\left(f^{\prime }\tau\right)-\frac{1+\left|
\tau\right| ^2}{2\left| f^{\prime }\right| }{{\rm{Re}}}\left((f^{\prime })^2\bar{f}
^{\prime \prime }\right),
\end{eqnarray*}
and similarly,
$$
\alpha _2 =-\left| f^{\prime }\right| {{\rm{Im}}}\left(f^{\prime }\tau\right)-\frac{
1+\left| \tau\right| ^2}{2\left| f^{\prime }\right| }{{\rm{Im}}}\left((f^{\prime })^2
\bar{f}^{\prime \prime }\right).
$$
On the other hand, 
\begin{eqnarray*}
\left| f^{\prime }\right| ^2+\left| \nabla R\right| ^2 & = &\left| f^{\prime
}\right| ^2+\left(u\left| f^{\prime }\right| +\frac{1+\left| \tau\right| ^2}{2\left|
f^{\prime }\right| }{{\rm{Re}}}\left(f^{\prime }\bar{f}^{\prime \prime
}\right)^2+v\left| f^{\prime }\right| +\frac{1+\left| \tau\right| ^2}{2\left|
f^{\prime }\right| }{{\rm{Im}}}\left(f^{\prime }\bar{f}^{\prime \prime }\right)\right)^2 \\
& = & 
\left(\left| \tau\right|
^2+1\right)\left(\left| f^{\prime }\right| ^2+\hbox{Re}\left( f^{\prime }\bar{f}
^{\prime \prime }\tau\right) +\frac{\left| f^{\prime \prime }\right| ^2(\left|
\tau\right| ^2+1)}4\right).
\end{eqnarray*}
Assembling together all the pieces of our computation we  can finally write down  the coordinates for
the sought CMC-1 surface as described in the theorem.
\end{proof}

\begin{remark}
Bianchi's construction is local in the sense that he assumes that the  projections defining $f$ are both bijective. But as the ruled example in Section \ref{examp} suggests, we may imagine a global construction by considering function elements on appropriate Riemann surfaces.
\end{remark}

\begin{remark}
Cal\`o also published formulae similar to (\ref{calo1})-(\ref{calo2}) corresponding to a rolling problem where the plane is replaced by the sphere (see \cite{b}). A variant of the argument above applied to this case then furnishes explicit formulae for CMC-1 surfaces in the Poincar\'e model. 

\end{remark}

\begin{remark}
Retracing through the above construction, it is not hard to check that $f$ is the hyperbolic Gauss map $G$ of the corresponding CMC-1 surface. More precisely, we have $f=G\circ X$ where $X=(x,y,z)$ is given by (\ref{bicalo}).
\end{remark}

We finish this section by briefly indicating how Small's result relates to the one presented here. Recall from (\ref{sma}) the local representation for the Gauss transform $\Gamma_{\Sc}$ of a null curve $\Sc\subset\sl$
in terms of a pair $(f,g)$ of holomorphic functions. According to (\ref{small}), one has to apply dualization to this in order to recover $\Sigma$. After doing this, one obtains that $\Sc$ is given by the map $\omega:M\to \sl$, 
\begin{equation}\label{ome}
\omega=\left(
\begin{array}{cc}
\alpha & \beta \\
\gamma & \delta
\end{array}
\right)=
\left( 
\begin{array}{cc}
(f')^{1/2}-\frac{1}{2}f(f')^{-3/2}f'' &
f\left((f')^{-1/2}+\frac{1}{2}g(f')^{-3/2}f''\right)-g(f')^{1/2} \\
-\frac{1}{2}(f')^{-3/2}f'' & 
(f')^{-1/2}+\frac{1}{2}g(f')^{-3/2}f''
\end{array}
\right),
\end{equation}
where $f'=df/dg$ and $f''=d^2f/dg^2$. The expression for the CMC-1 surface in terms of the hermitian model is, according to Bryant, given by $\o\overline{\o}^t$, but in order to compare this with (\ref{bicalo}) one has to perform the  transformation to the upper half-space model. In terms of the entries of $\o$, this  is given by 
$$
x+iy=\frac{\alpha\overline{\gamma}+\beta\overline{\delta}}{|\gamma|^2
+|\delta|^2},\,\,\,\,z=\frac{1}{|\gamma|^2
+|\delta|^2}.
$$
We now take $g=\tau$ and Small's $f$ to be our $f$. A straightforward computation yields the equivalence between the methods.

\section{Some examples and final remarks}
\label{examp}

We illustrate the method by retrieving two well-known examples. First, if we take 
$
f(\tau )=\tau ^2,
$
substitution in (\ref{bicalo}) after writing $\tau =re^{i\theta }$ yields 
\begin{eqnarray*}
x &=&-r^2\left( \cos 2\theta \right) \frac{5r^2+3}{7r^2+1}, \\
y &=&-r^2\left( \sin 2\theta \right) \frac{5r^2+3}{7r^2+1}, \\
z &=&\frac{8r^3}{7r^2+1}.
\end{eqnarray*}
This is a {\it catenoid cousin}.
A rough picture is given in figure 4.1

\medskip
\centerline{\epsfxsize=13.5cm \epsfysize=7.56cm \epsfbox{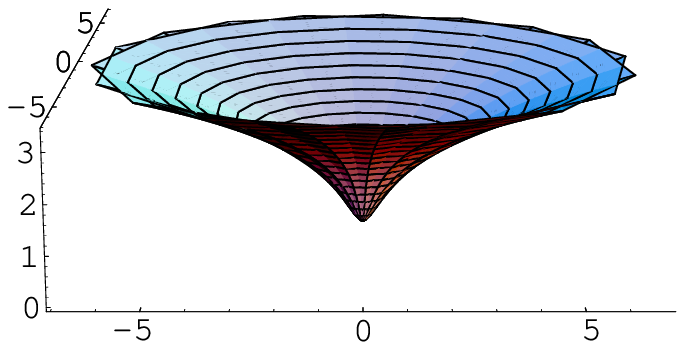}}   
\smallskip
\centerline{Fig. 4.1}
\medskip

Now let 
$
f(\tau )=\ln \tau.
$
Again, substitution in (\ref{bicalo}) yields 
\[
\mathbf{\xi }_{-}=(\ln r-2\frac{(r-r^{-1})}{(r+r^{-1})},\theta ,\frac
4{(r+r^{-1})}).
\]
Or, writing $r=e^s$,
\[
\mathbf{\xi }_{-}=(s-2\tanh s,\theta ,\frac 2{\cosh s}).
\]
This is a ruled example as figure 4.2  makes it clear.

\medskip
\centerline{\epsfxsize=13.5cm \epsfysize=7.56cm \epsfbox{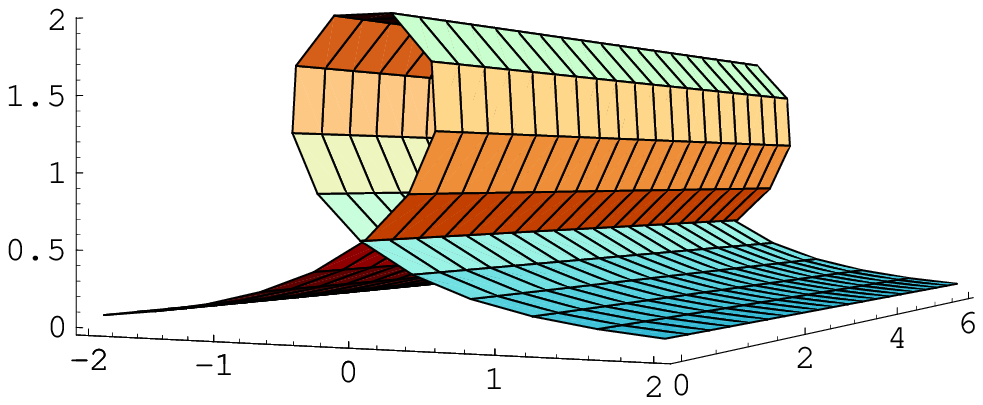}}   
\smallskip
\centerline{Fig. 4.2}
\medskip

\begin{remark}
Strictly speaking we used the Riemann surface associated to the $\log$ function in this last example, showing that the construction might work in more general cases.
\end{remark}

One  would like to add a few questions:

\begin{enumerate}
\item  Can one understand the interplay between the Euclidean problem 
(rolling of isometric surfaces) and the construction of CMC-1 surfaces in terms of some underlying structure? Can the relationship be thought of as a natural one from some other point of view?

\item Is there a way to view directly the congruence of spheres in a twistorial perspective (or in any other {complexified} way)?
\end{enumerate}

\end{document}